\documentclass{article}
\usepackage{fullpage}
\usepackage{fancyhdr}
\usepackage{amsmath, amssymb, latexsym, amscd, amsthm,amsfonts,amstext}
\usepackage{subfigure}
\usepackage{enumerate}
\usepackage{xcolor}
\usepackage{textcomp}
\usepackage{mathtools}

\def\div{\operatorname{div}}
\usepackage{graphicx}
\newtheorem{theorem}{Theorem}[section]
\newtheorem{lemma}[theorem]{Lemma}

\numberwithin{equation}{section}
\newtheorem{definition}[theorem]{Definition}

\numberwithin{equation}{section}

\usepackage[mathscr]{eucal}
\usepackage{graphicx}
\usepackage{subfigure}
\usepackage{cite}
\usepackage{url}
\usepackage{enumerate}
\usepackage{xcolor}

\def\div{\operatorname{div}}
\allowdisplaybreaks
 
\title{Location and size estimation of small rigid bodies using elastic far-fields}
\author{
Fadhel  Al-Musallam
\thanks{Department of Mathematics, Kuwait University, P.O. Box 13060, Safat, Kuwait. (Email: musallam@sci.kuniv.edu.kw).}
\and
 Durga Prasad Challa
\thanks{Department of mathematics, 
Inha university,  Incheon 402-751, S. Korea.
(Email: durga.challa@inha.ac.kr).
}
\and  Mourad Sini \footnotemark[1]
\thanks{RICAM, Austrian Academy of Sciences,
Altenbergerstrasse 69, A-4040, Linz, Austria.
(Email: mourad.sini@oeaw.ac.at, msini@sci.kuniv.edu.kw).
}
}

\begin{document}
\graphicspath{{Smallscattersreport_Figures_50incobs/}}
\maketitle
\begin{abstract}
 We are concerned with the linearized, isotropic and homogeneous elastic scattering problem by (possibly many) small rigid 
 obstacles of arbitrary Lipschitz regular shapes in 3D. 
 Based on the Foldy-Lax approximation, valid under a sufficient condition on the number of the obstacles, the size and 
the minimum distance between them, we show that any of the two body waves, namely the pressure waves P or the shear waves S,
 is enough for solving the inverse problem of detecting these scatterers and estimating their sizes.  Further, it is also shown that the shear-horizontal part SH or the shear vertical 
 part SV of the shear waves S are also enough for the location detection and the size estimation. Under some extra assumption on the scatterers, as the convexity assumption,
 we derive finer size estimates as the radius of the largest ball contained in each scatterer and the one of the smallest ball containing it. 
 The two estimates measure, respectively, the thickness and length of each obstacle. 
\end{abstract}


\section{Introduction and statement of the results}\label{Introduction-smallac-sdlp}

 Let $B_1, B_2,\dots, B_M$ be $M$ open, bounded and simply connected sets in $\mathbb{R}^3$ with Lipschitz boundaries,
 containing the origin. 
We assume that their sizes and Lipschitz constants are uniformly bounded.  
We set $D_m:=\epsilon B_m+z_m$ to be the small bodies characterized by the parameter 
$\epsilon>0$ and the locations $z_m\in \mathbb{R}^3$, $m=1,\dots,M$. 
\par 
 Assume that the Lam\'e coefficients $\lambda$ and $\mu$ are constants satisfying $ \mu > 0 \mbox{ and } 3\lambda+2\mu >0$.   
Let $U^{i}$ be a solution of the Navier equation $(\Delta^e + \omega^{2})U^{i}=0 \mbox{ in } \mathbb{R}^{3}$, $\Delta^{e}:=(\mu\Delta+(\lambda+\mu)\nabla \div)$. 
We denote by  $U^{s}$ the elastic field scattered by the $M$ small bodies $D_m\subset \mathbb{R}^{3}$ due to 
the incident field $U^{i}$. We restrict ourselves to the scattering by rigid bodies. Hence the total field $U^{t}:=U^{i}+U^{s}$ 
satisfies the following exterior Dirichlet problem of the elastic waves
\begin{equation}
(\Delta^e + \omega^{2})U^{t}=0 \mbox{ in }\mathbb{R}^{3}\backslash \left(\mathop{\cup}_{m=1}^M \bar{D}_m\right),\label{elaimpoenetrable}
\end{equation}
\begin{equation}
U^{t}|_{\partial D_m}=0,\, 1\leq m \leq M\label{elagoverningsupport}  
\end{equation}
with the Kupradze radiation conditions (K.R.C)
\begin{equation}\label{radiationcela}
\lim_{|x|\rightarrow\infty}|x|^{\frac{d-1}{2}}(\frac{\partial U_{p}}{\partial|x|}-i \kappa_{p^\omega}U_{p})=0,  \mbox{  and  } 
\lim_{|x|\rightarrow\infty}|x|^{\frac{d-1}{2}}(\frac{\partial U_{s}}{\partial|x|}-i \kappa_{s^\omega}U_{s})=0, 
\end{equation}
where the two limits are uniform in all the directions $\hat{x}:=\frac{x}{|x|}\in \mathbb{S}^{2}$ and $\mathbb{S}^{2}$ is the unit sphere. 
Also, we denote $U_{p}:=-\kappa_{p^\omega}^{-2}\nabla (\nabla\cdot U^{s})$ to be the longitudinal 
(or the pressure or P) part of the field $u$ and $U_{s}:=\kappa_{s^\omega}^{-2}\nabla\times(\nabla\times U^{s})$ to be the transversal (or the shear or S) part of 
the field $U^{s}$ corresponding to the Helmholtz decomposition $U^{s}=U_{p}+U_{s}$. The constants $\kappa_{p^\omega}:=\frac{\omega}{c_p}$ and  
$\kappa_{s^\omega}:=\frac{\omega}{c_s}$ are known as the longitudinal and transversal wavenumbers, $c_p:=\sqrt{\lambda+2\mu}$ and $c_s:=\sqrt{\mu}$ are the corresponding phase velocities, respectively and $\omega$ is the frequency.
\par 
The scattering problem (\ref{elaimpoenetrable}-\ref{radiationcela}) is well posed in the H\"{o}lder or Sobolev spaces, 
see \cite{A-B-G-K-L-W:2014, D-K:2000, K-G-B-B:1979} for instance, and the scattered field $U^s$ has the following asymptotic expansion:

\begin{equation}\label{Lamesystemtotalfieldasymptoticsmall}
 U^s(x) := \frac{e^{i\kappa_{p^\omega}|x|}}{|x|}U^{\infty}_{p}(\hat{x}) + 
\frac{e^{i\kappa_{s^\omega}|x|}}{|x|}U^{\infty}_{s}(\hat{x}) + O(\frac{1}{|x|^{2}}),~|x|\rightarrow \infty
\end{equation}
uniformly in all directions $\hat{x}\in \mathbb{S}^{d-1}$. The longitudinal part of the far-field, i.e. $U^{\infty}_{p}(\hat{x})$ is normal to $\mathbb{S}^{2}$ 
while the transversal part $U^{\infty}_{s}(\hat{x})$ is tangential to $\mathbb{S}^{2}$.  
As usual in scattering problems we use plane incident waves in this work. For the Lam\'e system, the full plane incident wave is of the form
$U^{i}(x,\theta):=\alpha\theta\,e^{i\kappa_{p^\omega}\theta\cdot x}+\beta\theta^\bot\,e^{i\kappa_{s^\omega}\theta\cdot x}$, 
where $\theta^{\bot}$ is any direction in $\mathbb{S}^{2}$ perpendicular to the incident direction $\theta \in \mathbb{S}^2$, $\alpha,\beta$ 
are arbitrary constants.  In particular, the pressure and shear incident waves are given as follows;  

\begin{equation}\label{elaincidentwaves}
 U^{i,p}(x,\theta) := \theta e^{i\kappa_{p^\omega}\theta\cdot x}\mbox{ and } U^{i,s}(x,\theta) := \theta^{\bot} e^{i\kappa_{s^\omega}\theta\cdot x}.
\end{equation}
 Pressure incident waves propagate in the direction of $\theta$, whereas shear 
incident waves propagate in the direction of $\theta^{\bot}$.  In the two dimensional case, the shear waves have only one direction.  But in the three dimensional case, 
they have two orthogonal components called vertical and horizontal shear directions denoted by $\theta^{\bot_{v}}$ and  $\theta^{\bot_{h}}$ respectively.  
So, $\theta^{\bot}=\theta^{\underline{\bot}}\slash|\theta^{\underline{\bot}}|$ with $\theta^{\underline{\bot}}:=\alpha \theta^{\bot_{h}}+ \beta \theta^{\bot_{v}}$ for arbitrary constants $\alpha$ and $\beta$.  
To give the explicit forms of $\theta^{\bot_{h}}$ and $\theta^{\bot_{v}}$, we recall the Euclidean basis $\{e_1,e_2,e_3\}$ where $e_1:=(1,0,0)^T,e_2:=(0,1,0)^T$ and $e_3:=(0,0,1)^T$, 
write $\theta:=(\theta_x,\theta_y,\theta_z)^T$ and set  $r^{2}:=\theta_x^{2}+\theta_y^{2}$.  
Let $\mathcal{R}_3=\mathcal{R}_3(\theta)$ be the rotation map transforming $\theta$ to $e_{3}$. Then in the basis $\{e_1,e_2,e_3\}, \mathcal{R}_3$ is given by the matrix

\begin{equation}\label{rotation}
 \mathcal{R}_3= \frac{1}{r^{2}}\left[\begin{array}{ccc}
    \theta^2_2 + \theta^2_1 \theta_z &-\theta_x \theta_y(1-\theta_z)    &  -\theta_x r^{2} \\
 -\theta_x \theta_y(1-\theta_z)  & \theta^2_1 + \theta^2_2 \theta_z  & -\theta_y r^{2}   \\
\theta_x r^{2}  & \theta_y r^{2}  & ~~\theta_z r^{2}
   \end{array}\right].
\end{equation}
It satisfies  $\mathcal{R}^T_3 \mathcal{R}_3=I$ and  $\mathcal{R}_3\theta=e_{3}$.  
Correspondingly, we write $\theta^{\bot_{h}}:=\mathcal{R}_3^Te_{1}$ and $\theta^{\bot_{v}}:=\mathcal{R}_3^Te_{2}$.
 These two  directions represent the horizontal and the vertical directions of the shear wave and they are given by
\begin{equation}\label{defthetavh}
 \theta^{\bot_{h}}=\frac{1}{r^{2}}(\theta_y^{2}+\theta_x^{2}\theta_z,\theta_x \theta_y(\theta_z-1),-r^{2}\theta_x)^\top,\hspace{.25cm}
  \theta^{\bot_{v}}=\frac{1}{r^{2}}(\theta_x \theta_y(\theta_z-1),\theta_x^{2}+\theta_y^{2}\theta_z,-r^{2}\theta_y)^\top.
\end{equation}
The functions
$U^{\infty}_p(\hat{x}, \theta):=U^{\infty}_p(\hat{x})$ and $U^{\infty}_s(\hat{x}, \theta):=U^{\infty}_s(\hat{x})$ for $(\hat{x}, \theta)\in \mathbb{S}^{2} \times\mathbb{S}^{2}$ 
are called the P-part and the S-part of the far-field pattern respectively.
\newline
\begin{definition} 
\label{Def1}
We define 
\begin{enumerate}
 \item $ a:=\max\limits_{1\leq m\leq M } diam (D_m) ~~\big[=\epsilon \max\limits_{1\leq m\leq M } diam (B_m)\big],$
 \item $d:=\min\limits_{\substack{m\neq j\\1\leq m,j\leq M }} d_{mj},$ where $d_{mj}:=dist(D_m, D_j)$.
\item $\omega_{\max}$ as the upper bound of the used frequencies, i.e. $\omega\in[0,\,\omega_{\max}]$.
\item $\Omega$ to be a bounded domain in $\mathbb{R}^3$ containing the small bodies $D_m,\,m=1,\dots,M$.
\end{enumerate}
\end{definition}

Our goal in this work is to justify the following results.
\begin{theorem}\label{Main-Result}
 The matrix $(U^{\infty}_p(\hat{x}_j, \theta_l))^{N}_{j,l =1}$ (or $(U^{\infty}_s(\hat{x}_j, \theta_l))^{N}_{j,l =1}$), for $N$ large enough, corresponding to one
 of the incident waves in (\ref{elaincidentwaves}) is enough to localize the centers $z_j$ and estimate the sizes of the obstacles $D_j$, $j=1, ..., M$.
 
 Further, the matrix $(U^{\infty}_{SV}(\hat{x}_j, \theta_l))^{N}_{j,l =1}$ (or $(U^{\infty}_{SH}(\hat{x}_j, \theta_l))^{N}_{j,l =1}$) is also enough to localize the obstacles and 
 estimate their sizes. Here $U^{\infty}_{SV}(\cdot, \cdot)$ and $U^{\infty}_{SH}(\cdot, \cdot)$ are respectively the Shear-Horizontal and the Shear-Vertical parts of the shear
 parts of the far-fields.
\end{theorem}

The approach we use to justify these results is based on two steps.
\begin{enumerate}
 \item In the first step, we derive the asymptotic expansion of the far-fields in terms of the three parameters modeling the collection of scatterers, namely $M$, 
 $a$ and $d$. This is sometimes called the Foldy-Lax approximation.
 
 \item We use the dominant term of this approximation coupled with the so-called MUSIC algorithm to detect the locations of the obstacles. 
 
\end{enumerate}
This approach is
 known for a decade, see \cite{A-K:2007} for instance. Our contribution to this approach is twofold corresponding to the two steps mentioned above. Regarding the first step, we
 provide the asymptotic expansion in terms of the three parameters $M$, $a$ and $d$, while in the previous literature the two parameters $M$ and $d$ are assumed to be 
 fixed, \cite{A-C-I:SIAM2008, A-B-G-K-L-W:2014}. This expansion is justified for the Lam\'e model, under consideration here, in our previous work \cite{DPC-SM13-2}. Regarding the second
 step, which is the object of this paper, we apply the MUSIC algorithm, see \cite{A-K:2007, DPC-SM:InvPrbm2012}, to the P-parts (respectively the S-parts) of the elastic far-fields to localize
 the centers of the scatterers. Further, we extract the elastic capacitances of the obstacles from these data. Finally, from these capacitances we derive lower and upper 
 estimates of the scaled perimeter of the scatterers, see Theorem \ref{general-obstacles}. If in addition the obstacles are convex, then we derive an upper estimate of the largest ball contained in each obstacle and
 a lower bound of the smallest ball containing it, see Theorem \ref{convex-obstacles}. The two estimates measure, respectively, the thickness and length of each obstacle. 
 It seems to us that these two estimates are new in the literature.
 \newline
 
 Let us also emphasize that our results mean that any of the two body waves (pressure or shear waves) is enough to localize and 
 estimate the sizes of the scatterers.
 \newline
 
The rest of paper is organized as follows. In section 2, we recall, from \cite{DPC-SM13-2}, the Foldy-Lax approximation of the elastic fields. In section 3, we use these
approximations to justify Theorem \ref{Main-Result}.

\section{Forward Problem} 
\subsection{The asymptotic expansion of the far-fields}
The forward problem is to compute the P-part, $U^\infty_p(\hat{x},\theta)$, and the S-part, $U^\infty_s(\hat{x},\theta)$, 
of the far-field pattern associated with the Lam\'e system (\ref{elaimpoenetrable}-\ref{radiationcela}) for various incident and the observational directions. 
The main result is the following theorem, see \cite[Theorem 1.2]{DPC-SM13-2}, which 
justifies the Foldy-Lax approximation, in order to represent the scattering by small scatterers taking into account the three parameters $M$, $a$ and $d$.

\begin{theorem}\label{Maintheorem-ela-small-sing}
 There exist two positive constants $a_0$ and $c_0$ depending only on the size of $\Omega$, the 
Lipschitz character of $B_m,m=1,\dots,M$, $d_{\max}$ and $\omega_{\max}$ such that
if 
\begin{equation}\label{conditions-elasma}
a \leq a_0 ~~ \mbox{and} ~~ \sqrt{M-1}\frac{a}{d}\leq c_0
\end{equation} 
then the P-part, $U^\infty_p(\hat{x},\theta)$, and the S-part, $U^\infty_s(\hat{x},\theta)$, of the far-field pattern have the following asymptotic expressions
 \begin{eqnarray}
  U^\infty_p(\hat{x},\theta)&=&\hspace{-.1cm}\frac{1}{4\pi\,c_p^{2}}(\hat{x}\otimes\hat{x})\hspace{-.1cm}\left[\sum_{m=1}^{M}e^{-i\frac{\omega}{c_p}\hat{x}\cdot z_{m}}Q_m\right.\left.+O\left(M a^2+M(M-1)\frac{a^3}{d^2}
+M(M-1)^2\frac{a^4}{d^3}\right)\right], \label{x oustdie1 D_m farmainp}\\
 U^\infty_s(\hat{x},\theta)&=&\hspace{-.1cm} \frac{1}{4\pi\,c_s^{2}}(I- \hat{x}\otimes\hat{x})\hspace{-.1cm}\left[\sum_{m=1}^{M}e^{-i\frac{\omega}{c_s}\hat{x}\cdot\,z_m}Q_m\right.\left.+O\left(M a^2+M(M-1)\frac{a^3}{d^2}
+M(M-1)^2\frac{a^4}{d^3}\right)\right]. \label{x oustdie1 D_m farmains}
  \end{eqnarray}
uniformly in $\hat{x}$ and $\theta$ in $\mathbb{S}^2$. The constant appearing in the estimate $O(.)$ depends only on the size of $\Omega$, 
the Lipschitz character of the reference bodies, $a_0$, $c_0$ and $\omega_{max}$. The vector coefficients $Q_m$, $m=1,..., M,$ are the solutions of the following linear algebraic system
\begin{eqnarray}\label{fracqcfracmain}
 C_m^{-1}Q_m+\sum_{\substack{j=1 \\ j\neq m}}^{M} \Gamma^{\omega}(z_m,z_j)Q_j &=&-U^{i}(z_m, \theta),~~
\end{eqnarray}
for $ m=1,..., M,$ with $\Gamma^{\omega}$ denoting the Kupradze matrix of the fundamental solution to the Navier equation with frequency 
$\omega$, $C_m:=\int_{\partial D_m}\sigma_m(s)ds$ and $\sigma_{m}$ is 
the solution matrix of the integral equation of the first kind
\begin{eqnarray}\label{barqcimsurfacefrm1main}
\int_{\partial D_m}\Gamma^{0}(s_m,s)\sigma_{m} (s)ds&=&\rm \textbf{I},~ s_m\in \partial D_m,
\end{eqnarray}
with $\rm \textbf{I}$ the identity matrix of order 3. The algebraic system \eqref{fracqcfracmain} is invertible 
under the condition:
\begin{eqnarray}\label{invertibilityconditionsmainthm-ela}
\frac{a}{d}&\leq&c_1t^{-1}
\end{eqnarray}
with 
\begin{center}
$
t:=\left[\frac{1}{c_p^2}-2diam(\Omega)\frac{\omega}{c_s^3}\left(\frac{1-\left(\frac{1}{2}\kappa_{s^\omega}diam(\Omega)\right)^{N_\Omega}}{1-\left(\frac{1}{2}\kappa_{s^\omega}diam(\Omega)\right)}+\frac{1}{2^{N_{\Omega}-1}}\right)-diam(\Omega)\frac{\omega}{c_p^3}\left(\frac{1-\left(\frac{1}{2}\kappa_{p^\omega}diam(\Omega)\right)^{N_\Omega}}{1-\left(\frac{1}{2}\kappa_{p^\omega}diam(\Omega)\right)}+\frac{1}{2^{N_{\Omega}-1}}\right)\right],
$
\end{center}
which is assumed to be positive
and $N_{\Omega}:=[2diam(\Omega)\max\{\kappa_{s^\omega},\kappa_{p^\omega}\}e^2]$, where $[\cdot]$ denotes the integral part and $\ln e=1$. 
The constant $c_1$ depends only on the Lipschitz character of the reference bodies $B_m$, $m=1,\dots,M$.
\end{theorem}

We call the system \eqref{fracqcfracmain} the elastic Foldy-Lax algebraic system. The matrix $C_m:=\int_{\partial D_m}\sigma_m(s)ds$, 
where $\sigma_m$ solves (\ref{barqcimsurfacefrm1main}), is called
the elastic capacitance of the set $D_m$. One of the interests of the expansions (\ref{x oustdie1 D_m farmainp}) and (\ref{x oustdie1 D_m farmains}) is that
we can reduce the computation of the elastic fields due to small obstacles to solving an algebraic system (i.e. \eqref{fracqcfracmain}) and inverting a 
first kind integral equation (i.e. \eqref{barqcimsurfacefrm1main}). Another goal in deriving the expansion in terms of the three parameters is the 
quantification of the equivalent effective medium, 
 without homogeneization (i.e. with no periodicity assumption on the distribution of the scatterers), see \cite{B-M-D-M} for the acoustic model.

\subsection{The fundamental solution}\label{fdelsms}
The Kupradze matrix $\Gamma^\omega=(\Gamma^\omega_{ij})^3_{i,j=1}$ of the fundamental solution to the Navier equation is given by
\begin{eqnarray}\label{kupradzeten}
 \Gamma^\omega(x,y)=\frac{1}{\mu}\Phi_{\kappa_{s^\omega}}(x,y)\rm \textbf{I}+\frac{1}{\omega^2}\nabla_x\nabla_x^{\top}[\Phi_{\kappa_{s^\omega}}(x,y)-\Phi_{\kappa_{p^\omega}}(x,y)],
\end{eqnarray}
where $\Phi_{\kappa}(x,y)=\frac{\exp(i\kappa|x-y|)}{4\pi |x-y|}$ denotes the free space fundamental solution of the Helmholtz equation $(\Delta+\kappa^2)\,u=0$ in 
$\mathbb{R}^3$.
The asymptotic behavior of Kupradze tensor at infinity is given as follows
\begin{equation}\label{elafundatensorasymptotic}
 \Gamma^{\omega}(x,y)=\frac{1}{4\pi\,c_p^{2}}\hat{x}\otimes\hat{x} \frac{e^{i\kappa_{p^\omega}|x|}}{|x|}e^{-i\kappa_{p^\omega}\hat{x}\cdot\,y} + 
\frac{1}{4\pi\,c_s^{2}}(I- \hat{x}\otimes\hat{x}) \frac{e^{i\kappa_{s^\omega}|x|}}{|x|}e^{-i\kappa_{s^\omega}\hat{x}\cdot\,y}+O(|x|^{-2})
\end{equation}
with $\hat{x}=\frac{x}{|x|}\in\mathbb{S}^{2}$, see \cite{D-K:2000} for instance.

\section{Inverse problem}\label{Inverseproblem}
\subsection{Scalar far-field patterns}\label{Scalar far-field patterns}
  We define the scalar P-part, $U^\infty_p(\hat{x},\theta)$, and the scalar S-part, $U^\infty_s(\hat{x},\theta)$, of the far-field pattern of the problem 
  (\ref{elaimpoenetrable}-\ref{radiationcela}) respectively as
 \begin{eqnarray}
   \mathrm{U}^\infty_p(\hat{x},\theta)&:=&4\pi\,c_p^{2}\, \left(\hat{x}\cdot U^\infty_p(\hat{x},\theta)\right)\nonumber\\
&=&\sum_{m=1}^{M}\hat{x}e^{-i\frac{\omega}{c_p}\hat{x}\cdot z_{m}}Q_m+O\left(M a^2+M(M-1)\frac{a^3}{d^2}+M(M-1)^2\frac{a^4}{d^3}\right),  \label{x oustdie1 D_melaP1} \\
 \mathrm{U}^\infty_s(\hat{x},\theta)&:=&4\pi\,c_s^{2}\,\left(\hat{x}^{\bot}\cdot U^\infty_s(\hat{x},\theta)\right)\nonumber\\
&=& \sum_{m=1}^{M}\hat{x}^{\bot}e^{-i\frac{\omega}{c_s}\hat{x}\cdot\,z_m}Q_m+O\left(M a^2+M(M-1)\frac{a^3}{d^2}+M(M-1)^2\frac{a^4}{d^3}\right). \label{x oustdie1 D_melaS1}
\end{eqnarray}
From \eqref{x oustdie1 D_melaP1} and \eqref{x oustdie1 D_melaS1}, we can write the scalar P and the scalar S parts of the far-field pattern as
 \begin{eqnarray}
   \mathrm{U}^\infty_p(\hat{x},\theta)&=&\sum_{m=1}^{M}\hat{x}e^{-i\frac{\omega}{c_p}\hat{x}\cdot z_{m}}Q_m\label{x oustdie1 D_melaP2}\\
 \mathrm{U}^\infty_s(\hat{x},\theta)&=&\sum_{m=1}^{M}\hat{x}^{\bot}e^{-i\frac{\omega}{c_s}\hat{x}\cdot\,z_m}Q_m\label{x oustdie1 D_melaS2}
\end{eqnarray}
with the error of order $O\left(M a^2+M(M-1)\frac{a^3}{d^2}+M(M-1)^2\frac{a^4}{d^3}\right)$
and $Q_m$ can be obtained from the linear algebraic system \eqref{fracqcfracmain}.  
\subsection{The elastic Foldy-Lax algebraic system}\label{sec-algebraicsys-small}
We can rewrite the  algebraic system \eqref{fracqcfracmain},
  \begin{eqnarray}\label{fracqcfracela}
 C_m^{-1}Q_m &=&-U^{i}(z_m)-\sum_{\substack{j=1 \\ j\neq m}}^{M} \Gamma^{\omega}(z_m,z_j)C_j(C_j^{-1}Q_j),
  \end{eqnarray}
 for all $m=1,2,\dots,M$. 
It can be written in a compact form as
\begin{equation}\label{compacfrm1ela}
 \mathbf{B}Q=U^I,
\end{equation}
\noindent
where $Q,U^I \in \mathbb{C}^{3M\times 1}\mbox{ and } \mathbf{B}\in\mathbb{C}^{3M\times 3M}$ are defined as
\begin{eqnarray}
\mathbf{B}:=\left(\begin{array}{ccccc}
   -{C}_1^{-1} &-\Gamma^{\omega}(z_1,z_2)&-\Gamma^{\omega}(z_1,z_3)&\cdots&-\Gamma^{\omega}(z_1,z_M)\\
-\Gamma^{\omega}(z_2,z_1)&-{C}_2^{-1}&-\Gamma^{\omega}(z_2,z_3)&\cdots&-\Gamma^{\omega}(z_2,z_M)\\
 \cdots&\cdots&\cdots&\cdots&\cdots\\
-\Gamma^{\omega}(z_M,z_1)&-\Gamma^{\omega}(z_M,z_2)&\cdots&-\Gamma^{\omega}(z_M,z_{M-1}) &-{C}_M^{-1}
   \end{array}\right),\label{The-matrix-B}\\
\nonumber\\
 Q:=\left(\begin{array}{cccc}
    Q_1^\top & Q_2^\top & \ldots  & Q_M^\top
   \end{array}\right)^\top \text{ and } 
U^I:=\left(\begin{array}{cccc}
     U^i(z_1)^\top & U^i(z_2)^\top& \ldots &  U^i(z_M)^\top
   \end{array}\right)^\top.
\nonumber
\end{eqnarray}
The above linear algebraic system is solvable for the 3D vectors $Q_j,~1\leq j\leq M$, when the matrix $\mathbf{B}$ is invertible. 
The invertibility of $\mathbf{B}$ is discussed in \cite[Corollary 4.3]{DPC-SM13-2}. 
\par Let us denote the inverse of $\mathbf{B}$ by $\mathcal{B}$ and the corresponding $3\times3$ blocks of $\mathcal{B}$ by $\mathcal{B}_{mj}$, $m,j=1,\dots,M$. 
Then we can rewrite \eqref{x oustdie1 D_melaP2} and \eqref{x oustdie1 D_melaS2}, with the same error, as follows
\begin{eqnarray}
 \mathrm{U}^\infty_p(\hat{x},\theta)&=&\sum_{m=1}^{M}\sum_{j=1}^{M}e^{-i\frac{\omega}{c_p}\hat{x}\cdot z_{m}}\hat{x}^{\top}\mathcal{B}_{mj}U^{i}(z_j,\theta)\label{x oustdie1 D_melaP3}\\
\mathrm{U}^\infty_s(\hat{x},\theta)&=&\sum_{m=1}^{M}\sum_{j=1}^{M}e^{-i\frac{\omega}{c_s}\hat{x}\cdot z_{m}}(\hat{x}^{\bot})^{\top}\mathcal{B}_{mj}U^{i}(z_j,\theta)\label{x oustdie1 D_melaS3}
 \end{eqnarray}
for a given incident direction $\theta$ and observation direction $\hat{x}$. 
From \eqref{x oustdie1 D_melaP3} and \eqref{x oustdie1 D_melaS3}, we can get the scalar P and the scalar S parts of 
the far-field patterns corresponding to plane incident P-wave $U^{i,p}(x,\theta)$ and S-wave $U^{i,s}(x,\theta)$, that we denote respectively by 
$\mathrm{U}^{\infty,p}_p(\hat{x},\theta),\, \mathrm{U}^{\infty,p}_s(\hat{x},\theta),\,\mathrm{U}^{\infty,s}_p(\hat{x},\theta),\, 
\mathrm{U}^{\infty,s}_s(\hat{x},\theta)$ as below
\begin{eqnarray}
 \mathrm{U}^{\infty,p}_p(\hat{x},\theta)&=&\sum_{m=1}^{M}\sum_{j=1}^{M}e^{-i\frac{\omega}{c_p}\hat{x}\cdot z_{m}}\hat{x}^{\top}\mathcal{B}_{mj}\theta\,e^{i\frac{\omega}{c_p}\theta\cdot z_{j}},\label{x oustdie1 D_melaPP}\\
\mathrm{U}^{\infty,p}_s(\hat{x},\theta)&=&\sum_{m=1}^{M}\sum_{j=1}^{M}e^{-i\frac{\omega}{c_s}\hat{x}\cdot z_{m}}(\hat{x}^{\bot})^{\top}\mathcal{B}_{mj}\theta\,e^{i\frac{\omega}{c_p}\theta\cdot z_{j}},\label{x oustdie1 D_melaPS}\\
 \mathrm{U}^{\infty,s}_p(\hat{x},\theta)&=&\sum_{m=1}^{M}\sum_{j=1}^{M}e^{-i\frac{\omega}{c_p}\hat{x}\cdot z_{m}}\hat{x}^{\top}\mathcal{B}_{mj}\theta^{\bot}\,e^{i\frac{\omega}{c_s}\theta\cdot z_{j}},\label{x oustdie1 D_melaSP}\\
\mathrm{U}^{\infty,s}_s(\hat{x},\theta)&=&\sum_{m=1}^{M}\sum_{j=1}^{M}e^{-i\frac{\omega}{c_s}\hat{x}\cdot z_{m}}(\hat{x}^{\bot})^{\top}\mathcal{B}_{mj}\theta^{\bot}\,e^{i\frac{\omega}{c_s}\theta\cdot z_{j}}.\label{x oustdie1 D_melaSS}
 \end{eqnarray}
All the far-field patterns (\ref{x oustdie1 D_melaPP}-\ref{x oustdie1 D_melaSS}) are valid with the same error which is equal to the error in (\ref{x oustdie1 D_melaP2}-\ref{x oustdie1 D_melaS2}). 
Now onwards, let $U^{\infty}(\hat{x},\theta)$ represents any one of the scattered fields mentioned above.\\ \ ~ \

\subsection{Localization of $D_{m}$'s via the MUSIC algorithm}\label{LocalisationviaMUSIC-smallac-sdlp-ela}

The MUSIC algorithm is a method to determine the locations $z_{m},m=1,2,\dots,M$, of the scatterers $D_m,m=1,2,\dots,M$ from the measured far-field pattern $U^{\infty}(\hat{x},\theta)$ for 
a finite set of incidence and observation directions, i.e. $\hat{x},\theta \in \{ \theta_{j},j=1,\dots,N\}\subset\mathbb{S}^{2}$.  
We refer the reader to the monograph \cite{K-G:2008} for more information about this algorithm.  
We follow the way in \cite{DPC-SM:InvPrbm2012} which is based on the presentation in \cite{K-G:2008}.  
\subsubsection{The factorization of the response matrix}\label{The factorization of the response matrix}

We assume that the number of scatterers is not larger than the number of incident and observation directions, precisely $N\geq 3M$.  
We define the response matrix $F\in\mathbb{C}^{N\times N}$ by
\begin{equation}\label{respomatdefela}
 F_{jl} := U^{\infty}(\theta_{j},\theta_{l}).
\end{equation}
The P-part and the S-part of the response matrix $F$ by ${F}_{p}$ and ${F}_{s}$ respectively. From (\ref{x oustdie1 D_melaP3}-\ref{x oustdie1 D_melaS3}) and \eqref{respomatdefela}, we can write
\begin{eqnarray}
 (F_{p})_{jl}&:=&U^\infty_p(\theta_{j},\theta_{l})\,=\,\sum_{m=1}^{M}\sum_{j=1}^{M}e^{-i\frac{\omega}{c_p}\theta_{j}\cdot z_{m}}\theta_{j}^{\top}\mathcal{B}_{mj}U^{i}(z_j,\theta_{l})\nonumber\\
&\hspace{-2.6cm}=&\hspace{-1.75cm}\left[\theta_{j}^{\top}e^{-i\frac{\omega}{c_p}\theta_j\cdot z_1},\,\theta_{j}^{\top}e^{-i\frac{\omega}{c_p}\theta_j\cdot z_2},\,\cdots,\,\theta_{j}^{\top}e^{-i\frac{\omega}{c_p}\theta_j\cdot z_M}\right]\hspace{-.1cm}\mathcal{B}\hspace{-.05cm}\left[(U^{i}(z_j,\theta_{l}))^{\top},\,(U^{i}(z_2,\theta_{l}))^{\top},\,\cdots,\,(U^{i}(z_m,\theta_{l}))^{\top}\right]^{\top}\hspace{-.3cm},\nonumber\label{x oustdie2 D_melaP3}\\
(F_{s})_{jl}&:=&U^\infty_s(\theta_{j},\theta_{l})\,=\,\sum_{m=1}^{M}\sum_{j=1}^{M}e^{-i\frac{\omega}{c_s}\theta_{j}\cdot z_{m}}(\theta_{j}^{\bot})^{\top}\mathcal{B}_{mj}U^{i}(z_j,\theta_{l})\nonumber\\
&\hspace{-2.6cm}=&\hspace{-1.75cm}\left[(\theta_{j}^{\bot})^{\top}e^{-i\frac{\omega}{c_s}\theta_j\cdot z_1},\,(\theta_{j}^{\bot})^{\top}e^{-i\frac{\omega}{c_s}\theta_j\cdot z_2},\,\cdots,\,(\theta_{j}^{\bot})e^{-i\frac{\omega}{c_s}\theta_j\cdot z_M}\right]\hspace{-.1cm}\mathcal{B}\hspace{-.05cm}\left[(U^{i}(z_j,\theta_{l}))^{\top},\,(U^{i}(z_2,\theta_{l}))^{\top},\,\cdots,\,(U^{i}(z_m,\theta_{l}))^{\top}\right]^{\top}\hspace{-.3cm}.\nonumber\label{x oustdie2 D_melaS3}
 \end{eqnarray}
for all $j,l=1,\dots,N$.
In PP, PS, SS and SP scatterings, denote the response matrix $F$ by ${F}^{p}_{p}$, ${F}^{p}_{s}$, ${F}^{s}_{s}$ and ${F}^{s}_{p}$ respectively and these can be factorized as
\begin{align}
 F^{p}_{p}={H^{p}}^{*}\mathcal{B}H^{p},~
F^{p}_{s}={H^{s}}^{*}\mathcal{B}H^{p},~
 F^{s}_{s}={H^{s}}^{*}\mathcal{B}H^{s},~
\mbox{ and }F^{s}_{p}={H^{p}}^{*}\mathcal{B}H^{s}.\label{elafactSPsmall}
\end{align}
Here, the matrices $H^{p}\in\mathbb{C}^{3M\times N}$ and $H^{s}\in\mathbb{C}^{3M\times N}$ are defined as,
\[
H^{p}:=\left(\begin{array}{cccc}
    \theta_{1}e^{i\frac{\omega}{c_p}\theta_{1}\cdot z_{1}}& \theta_{2}e^{i\frac{\omega}{c_p}\theta_{2}\cdot z_{1}}&\dots&\theta_{N}e^{i\frac{\omega}{c_p}\theta_{N}\cdot z_{1}}\\
\theta_{1}e^{i\frac{\omega}{c_p}\theta_{1}\cdot z_{2}}& \theta_{2}e^{i\frac{\omega}{c_p}\theta_{2}\cdot z_{2}}&\dots&\theta_{N}e^{i\frac{\omega}{c_p}\theta_{N}\cdot z_{2}}\\
\dots&\dots &\dots& \dots\\
\theta_{1}e^{i\frac{\omega}{c_p}\theta_{1}\cdot z_{M}}&\theta_{2}e^{i\frac{\omega}{c_p}\theta_{2}\cdot z_{M}} &\dots&\theta_{N}e^{i\frac{\omega}{c_p}\theta_{N}\cdot z_{M}}
\end{array}\right),
\]
and
\[
H^{s}:=\left(\begin{array}{cccc}
    \theta_{1}^{\bot}e^{i\frac{\omega}{c_s}\theta_{1}\cdot z_{1}}& \theta_{2}^{\bot}e^{i\frac{\omega}{c_s}\theta_{2}\cdot z_{1}}&\dots&\theta_{N}^{\bot}e^{i\frac{\omega}{c_s}\theta_{N}\cdot z_{1}}\\
\theta_{1}^{\bot}e^{i\frac{\omega}{c_s}\theta_{1}\cdot z_{2}}& \theta_{2}^{\bot}e^{i\frac{\omega}{c_s}\theta_{2}\cdot z_{2}}&\dots&\theta_{N}^{\bot}e^{i\frac{\omega}{c_s}\theta_{N}\cdot z_{2}}\\
\dots&\dots &\dots& \dots\\
\theta_{1}^{\bot}e^{i\frac{\omega}{c_s}\theta_{1}\cdot z_{M}}&\theta_{2}^{\bot}e^{i\frac{\omega}{c_s}\theta_{2}\cdot z_{M}} &\dots&\theta_{N}^{\bot}e^{i\frac{\omega}{c_s}\theta_{N}\cdot z_{M}}
 \end{array}\right).
\]
In order to determine the locations $z_{m}$, we consider a 3D-grid of sampling points $z\in\mathbb{R}^{3}$ in a region containing the scatterers $D_{1},D_{2},\dots,D_{M}$.
For each point $z$, we define the vectors $\phi_{z,p}^{j}$ and $\phi_{z,s}^{j}$ in $\mathbb{C}^{N}$ by
\begin{align}
 \phi^{j}_{z,p} &:=\left((\theta_{1}\cdotp e_{j})e^{-i\frac{\omega}{c_p}\theta_{1}\cdotp z},(\theta_{2}\cdotp e_{j}) e^{-i\frac{\omega}{c_p}\theta_{2}\cdotp z},\dots, 
(\theta_{N}\cdotp e_{j})e^{-i\frac{\omega}{c_p}\theta_{N}\cdotp z}\right)^{T},\label{elaphizMS2SPPsmall}\\
 \phi^{j}_{z,s} &:= \left((\theta_{1}^{\bot}\cdotp e_{j})e^{-i\frac{\omega}{c_s}\theta_{1}\cdotp z},(\theta_{2}^{\bot}\cdotp e_{j}) e^{-i\frac{\omega}{c_s}\theta_{2}\cdotp z},\dots, 
(\theta_{N}^{\bot}\cdotp e_{j})e^{-i\frac{\omega}{c_s}\theta_{N}\cdotp z}\right)^{T}, \forall j=1,2,3.\label{elaphizMS1SPsmall}
\end{align}
\noindent
\subsubsection{MUSIC characterization of the response matrix}\label{MUSICchar-smallela-sdlp} Recall that MUSIC is essentially based on characterizing the range of the response matrix (signal space), forming projections
onto its null (noise) spaces, and computing its singular value decomposition.  In other words,
the MUSIC algorithm is based on the property that the test vector $\phi_{z,r}^{j}$ is in the range $\mathcal{R}(F_{r})$ of $F_r$ if and only if $z$ is at one of locations of the scatterers, see \cite{DPC-SM:InvPrbm2012}. 
Here, $F_r:=F^{p}_{r}$ or $F_p:=F^{s}_{r}$ and $r\in\{p,s\}$.
\par It can be proved based on the non-singularity of the scattering matrix $\mathcal{B}$ in the factorizations (\ref{elafactSPsmall}) of $F^{r_1}_{r_2},\,r_1,r_2\in\{p,s\}$. 
Due to this, the standard linear algebraic argument yields that, if $N\geq 3M$ and the if the matrix $H^{r}$ has maximal rank $3M$, then the ranges $\mathcal{R}(H^{r^*})$ and $\mathcal{R}(F_r)$ coincide.
\par For sufficiently large number $N$ of incident and the observational directions by following the same lines as in 
\cite{K-G:2008,A-C-I:SIAM2008,DPC-SM:InvPrbm2012}, 
the maximal rank property of $H$ can be justified. In this case MUSIC algorithm is applicable for our response matrices $F_p^p,F_s^p,F_p^s$ and $F_s^s$. 
\par From the above discussion,  MUSIC characterization of the locations of the small scatterers in elastic exterior Drichlet problem can be written as the following

\begin{theorem}\label{elamusictheoremPPMSS-small}

For $N\geq 3M$ sufficiently large, we have
\begin{eqnarray}\label{elatheoremstatementMS1}
 z \in \{z_{1},...,z_{M}\} &\Longleftrightarrow~\phi^{j}_{z,t}\in \mathcal{R}({H^{t}}^{*}),\mbox{ for some } j=1,2,3 \mbox{ and for all } t\in\{p,s\}.
\end{eqnarray}
Furthermore, the ranges of ${H^{t}}^{*}$ and $F^{r}_{t}$ coincide and thus
\begin{equation}\label{elatheoremstatementMS2}
 z \in \{z_{1},...,z_{M}\} \Longleftrightarrow \phi^{j}_{z,t}\in \mathcal{R}(F^{r}_{t}) \Longleftrightarrow \mathcal{P}_t\phi^{j}_{z,t}=0, \mbox{ for some } j=1,2,3 \mbox{ and for all } r,t\in\{p,s\}
\end{equation}
where $\mathcal{P}_t : \mathbb{C}^{N}\rightarrow \mathcal{R}(F^{r}_{t})^{\bot}= \mathcal{N}({F^{r}_{t}}^{*})$ is the orthogonal projection onto the null space $\mathcal{N}({F^{r}_{t}}^{*})$ of ${F^{r}_{t}}^{*}$.
\end{theorem}
\ ~ \ \par
From Theorem \ref{elamusictheoremPPMSS-small}, the MUSIC algorithm holds for the response matrices corresponding to the PP, PS, SS and SP scatterings. To make the best use of the singular value decomposition
in SP and PS scatterings, we apply the MUSIC algorithm to $F^{s}_{p}{F^{s}_{p}}^{*}~ (\mbox{resp, }{F^{s}_{p}}^{*}F^{s}_{p})$ and 
${F^{p}_{s}}^{*}F^{p}_{s}~(\mbox{resp, }F^{p}_{s}{F^{p}_{s}}^{*})$ with the help of the test vectors $\phi^{j}_{z,p}~(\mbox{resp, } \phi^{j}_{z,s})$ respectively. 
  
\par As we are dealing with the 3D case, while dealing with S incident 
wave or S-part of the far-field pattern, it is enough to use one of its horizontal (S$^h$) or vertical (S$^v$) parts.  
Hence, it is enough to study the far-field pattern of any of the  PP, PS$^h$, PS$^v$, S$^h$S$^h$, S$^h$S$^v$, S$^v$S$^h$, S$^v$S$^v$, S$^h$P, 
S$^v$P elastic scatterings to locate the scatterers.
In other words, in three dimensional case, instead of using the full incident wave and the full far-field pattern, it is enough to study one combination 
of pressure (P), horizontal shear (S$^{h}$) or vertical shear (S$^{v}$) parts of the elastic incident wave and a corresponding part of the 
elastic far-field patterns, see \cite{DPC-SM:InvPrbm2012}.
\newline
Indeed, define the vectors $\phi^{j}_{z,s^h},\phi^{j}_{z,s^v}\in\mathbb{C}^{N}$ and the matrices $H^{s^h},H^{s^v}\in\mathbb{C}^{3M\times N}$ 
exactly as $\phi^{j}_{z,s}$ and $H^s$ replacing $\theta_i^{\bot}$ for $i=1,\dots, N$ by $\theta_i^{\bot_h}$ and
 $\theta_i^{\bot_v}$ respectively, see \eqref{defthetavh}. We denote the response matrices by  
$F^{p}_{s^h},~F^{s^h}_{p},~F^{p}_{s^v},~F^{s^v}_{p},~F^{s^h}_{s^h},~F^{s^h}_{s^v},~F^{s^v}_{s^h},~\mbox{and }F^{s^v}_{s^v}$ in 
the elastic PS$^h$, S$^h$P, PS$^v$, S$^v$P, S$^h$S$^h$, S$^h$S$^v$, S$^v$S$^h$, S$^v$S$^v$ scatterings respectively, then we can state the following theorem related to the MUSIC algorithm for sufficiently large number of incident and observation angles,
\begin{theorem}\label{elamusictheoremPPMShSv-small}
For $N\geq 3M$ sufficiently large, we have
\begin{eqnarray}\label{elatheoremstatementMShv1}
 z \in \{y_{1},\dots,y_{M}\} &\Longleftrightarrow~\phi^{j}_{z,t}\in \mathcal{R}({H^{t}}^{*}),\mbox{ for some } j=1,2,3 \mbox{ and for all } t\in\{p,s^h,s^v\}.
\end{eqnarray}
Furthermore, the ranges of ${H^{t}}^{*}$ and $F^{r}_{t}$ coincide and thus
\begin{equation}\label{elatheoremstatementMShv2}
 z \in \{y_{1},\dots,y_{M}\} \Longleftrightarrow \phi^{j}_{z,t}\in \mathcal{R}(F^{r}_{t}) \Longleftrightarrow 
\mathcal{P}_t\phi^{j}_{z,t}=0, \mbox{ for some } j=1,2,3 \mbox{ and for all } r,t\in\{p,s^h,s^v\}
\end{equation}
where $\mathcal{P}_t : \mathbb{C}^{N}\rightarrow \mathcal{R}(F^{r}_{t})^{\bot}= \mathcal{N}({F^{r}_{t}}^{*})$ is 
the orthogonal projection onto the null space $\mathcal{N}({F^{r}_{t}}^{*})$ of ${F^{r}_{t}}^{*}$.
\end{theorem}

The proof of the previous two theorems can be carried out in the same lines as in \cite{DPC-SM:InvPrbm2012}.

\subsection{Estimating the sizes}\label{Estimating the sizes}
\subsubsection{Recovering the capacitances}\label{Recovering the capacitances-smallela-sdlp}
 Once we locate the scatterers from the given far-field patterns using the MUSIC algorithm, we can recover the capacitances $C_{m}$ of $D_m$ 
from the factorization $F^{r}_{t}={H^{t}}^{*}\mathcal{B}H^{r}$ of $F^{r}_{t}\in\mathbb{C}^{N\times N},\,r,t\in\{p,s^h,s^v\}$.
Indeed,  we know that the matrix $H^{t}$ has maximal rank, see Theorem 3.1 and Theorem 3.2 of \cite{DPC-SM:InvPrbm2012}.  So, the matrix $H^{t}H^{t^*}\in\mathbb{C}^{3M\times 3M}$ is invertible.  
Let us denote its inverse by $I_{H^t}$. Once we locate the scatterers through finding the locations $z_{1},z_{2},\dots,z_{M}$ by using the MUSIC algorithm for the given far-field patterns, 
we can recover $I_{H^{t}}$ and hence the matrix $\mathcal{B}\in C^{3M\times 3M}$ given by $\mathcal{B}=I_{H^{t}}H^{t}F^{r}_{t}H^{r^*}I_{H^r}$, where $I_{H^{t}} H^{t}$ (resp, ${H^{r}}^{*}I_{H^{r}}$) is the pseudo inverse of ${H^{t}}^{*}$ (resp, ${H^{r}}$).
As we know the structure of $\mathbf{B}\in\mathbb{C}^{3M\times 3M}$, the inverse of $\mathcal{B}\in\mathbb{C}^{3M\times 3M}$, we can recover the capacitance matrices $C_{1},C_{2},\dots,C_{M}$ 
of the small scatterers $D_1,D_2,\dots,D_M$ from the diagonal blocks of $\mathbf{B}$, see (\ref{The-matrix-B}). 
From these capacitances, we can estimate the size of the obstacles as follows.

\subsubsection{Estimating the sizes of the obstacles from the capacitances}
Let us first start with the following lemma which compares the elastic and the acoustic capacitances \footnote{Recall that, for $m=1,\dots,M$, ${C}^a_m:=\int_{\partial D_m}\sigma_m(s)ds$ and $\sigma_{m}$ is 
the solution of the integral equation of the first kind $\int_{\partial D_m}\frac{\sigma_{m} (s)}{4\pi|t-s|}ds=1,~ t\in \partial D_m$, see \cite{D-M:MMS:2014, M-M:Book2013}.}, 
see \cite{DPC-SM13-2, M-M:Book2013}. 
\begin{lemma}\label{capacitance-eig-single} Let $\lambda^{min}_{eig_m}$ and $\lambda^{max}_{eig_m}$ be the minimal and maximal eigenvalues of the elastic capacitance matrices ${C}_m$, for $m=1,2,\dots,M$. Denote 
by ${C}^a_m$ the capacitance of each scatterer in the acoustic case, then we have the following estimate;
\begin{eqnarray}\label{lowerupperestforintgradub-m}
 \mu\,C^a_m\,\leq\,\lambda^{min}_{eig_m}\,\leq\,\lambda^{max}_{eig_m}\,\leq\,(\lambda+2\mu)\,C^a_m,\quad\text{for}\quad m=1,2,\dots,M.
\end{eqnarray}
\end{lemma}
Now, let us derive lower and upper bounds of the sizes  of the obstacles in terms of the acoustic capacitances. 

Assume that $D_j$'s are balls of radius $\rho_j$, 
and center $0$ for simplicity, then we know that $\int_{\{y:|y|=\rho_j\}}\frac{dS_y}{|x-y|}=4\pi\rho_j$, for $\vert x \vert = \rho_j$, as observed 
in \cite[formula (5.12)]{M-M:MathNach2010}. Hence $\sigma_j(s)=\rho^{-1}_j$ and then $C^{a}_{j}=\int_{\partial D_j}\rho^{-1}_jds=2\pi \rho_j$ from which we can estimate 
the radius $\rho_j$. Other geometries, as cylinders, for which one can estimate exactly the size, from the capacitance, are shown in chapter 4 of \cite{MRAMAG-book1}. 

For general geometries, we proceed as follows. First, we recall the following result, see \cite{DPC-SM13-2}.

\begin{lemma}\label{lemmadifssbQQbCCb1eladbl}
For every $1\leq j\leq M$, the capacitance $C^a_j$  is of the form
\begin{eqnarray}\label{asymptotCapeladbl}
{C}^a_j\,=\,{C}^a_{B_j}\epsilon
 \end{eqnarray}
where ${C}^a_{B_j}$ is the acoustic capacitance of $B_j$.
\end{lemma}
Now let us consider a single obstacle $D:=\epsilon\,B+z$. Since $C^a_B:=\int_{\partial B}\sigma(s)ds$ and $\int_{\partial B}\frac{\sigma(s)}{4\pi|t-s|}ds=1$, 
from the invertibility of the single layer potential $S: L^2(\partial B) \rightarrow H^1(\partial B)$, defined as $Sf(t)=\int_{\partial B}\frac{f(s)}{4\pi|t-s|}ds=1$, we deduce that
\begin{equation}\label{r-estimate}
C^a_B \leq \vert \partial B \vert^{\frac{1}{2}} \Vert \sigma \Vert_{L^2(\partial B)} \leq \vert \partial B \vert^{\frac{1}{2}} 
\Vert S^{-1}\Vert_{\mathcal{L}\left(H^1(\partial B), L^2(\partial B) \right)} \vert \partial B \vert^{\frac{1}{2}}=
\Vert S^{-1}\Vert_{\mathcal{L}\left(H^1(\partial B), L^2(\partial B) \right)} \vert \partial B \vert.
\end{equation}
On the other hand, we recall the following lower estimate, see Theorem 3.1 in \cite{MRAMAG-book1} for instance,
\begin{equation}\label{C-J}
C^a_B \geq \frac{4\pi \vert \partial B \vert^2}{J} 
\end{equation}
where $J:=\int_{\partial B}\int_{\partial B}\frac{1}{\vert s-t\vert}ds dt$. Remark that $J=4\pi\int_{\partial B} S(1)(s)ds$. Hence
$$ J\leq  4 \pi \vert \partial B\vert^{\frac{1}{2}} \Vert S\Vert_{\mathcal{L}\left(L^2(\partial B), H^1(\partial B)\right)}
\Vert 1\Vert_{H^1(\partial B)}\leq 4 \pi \Vert S\Vert_{\mathcal{L}\left(L^2(\partial B), H^1(\partial B)\right)} \vert
\partial B\vert
$$
and using (\ref{C-J}) we obtain the lower bound
\begin{equation}\label{l-estimate}
C^a_B\geq \Vert S\Vert^{-1}_{\mathcal{L}\left(L^2(\partial B), H^1(\partial B) \right)} \vert \partial B\vert.
\end{equation}
Finally combining (\ref{r-estimate}) and (\ref{l-estimate}), we derive the estimate
\begin{equation}\label{Capacitance-size-estimate}
 \Vert S\Vert^{-1}_{\mathcal{L}\left(L^2(\partial B), H^1(\partial B) \right)} \vert \partial B\vert \leq C^a_B \leq
\Vert S^{-1}\Vert_{\mathcal{L}\left(H^1(\partial B), L^2(\partial B) \right)} \vert \partial B \vert
\end{equation}

Using Lemma \ref{lemmadifssbQQbCCb1eladbl} and the relation $\vert \partial D_\epsilon \vert = \epsilon^2 \vert\partial B\vert $ we obtain
the following size estimation:
\begin{equation}\label{Capacitance-size-estimate-D-epsilon}
 \Vert S^{-1}\Vert_{\mathcal{L}\left(H^1(\partial B), L^2(\partial B) \right)}^{-1} \epsilon C^a_{\epsilon} \leq \vert \partial D_\epsilon \vert \leq 
 \Vert S\Vert_{\mathcal{L}\left(L^2(\partial B), H^1(\partial B) \right)} \epsilon C^a_{\epsilon}.
\end{equation}

Now, using Lemma \ref{capacitance-eig-single}, we derive the following lower and upper bounds of the sizes of the obstacles $D_m$

\begin{equation}\label{Capacitance-size-estimate-D-epsilon-Lame}
 \Vert S^{-1}\Vert_{\mathcal{L}\left(H^1(\partial B_m), L^2(\partial B_m) \right)}^{-1}  (\lambda+2\mu)^{-1}\lambda^{max}_{eig_m} \leq 
 \frac{\vert \partial D_m\vert}{\epsilon}  \leq 
 \Vert S\Vert_{\mathcal{L}\left(L^2(\partial B_m), H^1(\partial B_m) \right)} \mu^{-1} \lambda^{min}_{eig_m}.
\end{equation}
Observe that one can estimate the norms of the operators appearing in (\ref{Capacitance-size-estimate-D-epsilon-Lame}) in terms of (only) the Lipschitz character. 
We summarize this result in the following theorem
\begin{theorem}\label{general-obstacles}
 There exist two constants $c(Lip)$ and $C(Lip)$ depending only on the Lipschitz character of $B_1, ..., B_M$, such that
 \begin{equation}\label{size-estimation}
 c(Lip)  (\lambda+2\mu)^{-1}\lambda^{max}_{eig_m} \leq 
 \frac{\vert \partial D_m\vert}{\epsilon}  \leq C(Lip) \mu^{-1} \lambda^{min}_{eig_m}.
\end{equation}
 \end{theorem}
Precisely $C(Lip)$ and $c(Lip)$ are characterized respectively by 
$\Vert S\Vert_{\mathcal{L}\left(L^2(\partial B_m), H^1(\partial B_m) \right)} \leq C(Lip)$ and
$\Vert S^{-1}\Vert_{\mathcal{L}\left(H^1(\partial B_m), L^2(\partial B_m) \right)} \leq c^{-1}(Lip)$.
 We can use the estimate (\ref{Capacitance-size-estimate-D-epsilon-Lame}) to provide the lower and the upper estimates of the scaled 'size' 
of scatterers $\frac{\vert \partial D_m\vert}{\epsilon}$. Under some conditions of the reference obstacles, we can derive more explicit size estimates. Let us first define 
\begin{equation}
 \delta(x):=\min_{y \in \partial B}\vert x-y\vert,\;~~ x\in \mathbb{R}^3~~ \text{ and } R_i(B):=\sup_{x \in B}\delta(x).
\end{equation}
The quantity $R_i(B)$ is the radius of the largest ball contained in $B$. Now, we set 
\begin{equation}
 R_e(B):=\frac{1}{2}\max_{x, y \in \bar{B}}\vert x-y\vert.
\end{equation}
The quantity $R_e(B)$ is the radius of the smallest ball containing $B$.

By the Gauss theorem, we see that
\begin{equation}\label{Gauss}
 \vert B \vert =\frac{1}{3}\int_B div(x)\;dx = \frac{1}{3}\int_{\partial B}s\cdot n(s)ds
\end{equation}
hence $ \vert B \vert \leq \frac{1}{3} \vert \partial B\vert \max_{s \in \partial B}\vert s\vert \leq  \frac{2}{3} \vert \partial B\vert R_e(B)$
since $\max_{s \in \partial B}\vert s\vert \leq 2 R_e(B)$. Hence
\begin{equation}\label{lowerbound-partial B}
 \vert \partial B\vert \geq \frac{3}{2} \frac{\vert B\vert}{R_e(B)}
\end{equation}

To derive the upper bound for $\vert \partial B\vert$, we use the following argument borrowed from (\cite{Kovarik:2014}, section 4). If we assume that $B$ is convex, 
then $\delta(\cdot)$ is a concave function in $B$ and then, see \cite{B-V},
\begin{equation}\label{concavity}
 \nabla \delta(x)\cdot(y-x)\geq \delta(y)-\delta(x),\; ~~ x, y \in \bar{B}.
\end{equation}
But $\nabla \delta(x) =-\nu(x)$ for $x\in \partial B$, where $\nu$ is the external unit normal to $\partial B$.
Let us now assume, in addition to the convexity property, that $B(0, R_i(B)) \subset B$ which roughly means that the origin is the 'center' of $B$.
With this assumption, taking $y=0$ in (\ref{concavity}), we obtain $
s \cdot \nu(s) \geq \delta(0)\geq R_i(B),\,(\text{since}\, B(0, R_i(B)) \subset B$).
Replacing in (\ref{Gauss}), we obtain $
 \vert B\vert \geq \frac{1}{3} \vert \partial B\vert R_i(B)
$ hence

\begin{equation}\label{upperbound-partial B}
 \vert \partial B\vert \leq 3 \frac{\vert B\vert}{R_i(B)}.
\end{equation}

Replacing (\ref{lowerbound-partial B}) and (\ref{upperbound-partial B}) in (\ref{Capacitance-size-estimate}), we have
$$
\frac{3}{2}\Vert S\Vert^{-1}_{\mathcal{L}\left(L^2(\partial B), H^1(\partial B) \right)} \frac{\vert B\vert}{R_e(B)} \leq C^a_B \leq
3\Vert S^{-1}\Vert_{\mathcal{L}\left(H^1(\partial B), L^2(\partial B) \right)} \frac{\vert B \vert}{R_i(B)}.
$$
Using the double inequality $
\frac{4}{3}\pi R^3_i(B) \leq \vert B\vert \leq \frac{4}{3}\pi R^3_e(B)$
we obtain
\begin{equation}\label{key-observation-B}
2 \pi \Vert S\Vert^{-1}_{\mathcal{L}\left(L^2(\partial B), H^1(\partial B) \right)} \frac{R^2_i(B)}{R_e(B)} R_i(B) \leq C^a_B \leq
4 \pi\Vert S^{-1}\Vert_{\mathcal{L}\left(H^1(\partial B), L^2(\partial B) \right)} \frac{R^2_e(B)}{R_i(B)} R_e(B).
\end{equation}
Now, we apply these double estimates to $D$ instead of $B$, knowing that the two assumptions on $B$ are inherited by $D$, then we obtain
\begin{equation}\label{key-observation-D}
2 \pi \Vert S\Vert^{-1}_{\mathcal{L}\left(L^2(\partial D), H^1(\partial D) \right)} \frac{R^2_i(D)}{R_e(D)} R_i(D) \leq C^a_D \leq
4 \pi\Vert S^{-1}\Vert_{\mathcal{L}\left(H^1(\partial D), L^2(\partial D) \right)} \frac{R^2_e(D)}{R_i(D)} R_e(D).
\end{equation}

Observe that $\frac{R^2_i(D)}{R_e(D)}$ scales as $\Vert S\Vert_{\mathcal{L}\left(L^2(\partial D), H^1(\partial D) \right)}$ and $\frac{R^2_e(D)}{R_i(D)}$ scales as
$\Vert S^{-1}\Vert^{-1}_{\mathcal{L}\left(H^1(\partial B), L^2(\partial B) \right)}$ since obviously $R_i(D)=\epsilon R_i(B), \; R_e(D)=\epsilon R_e(B)$ and
we have $\Vert S\Vert_{\mathcal{L}\left(L^2(\partial D), H^1(\partial D) \right)}\leq\epsilon \Vert S\Vert_{\mathcal{L}\left(L^2(\partial B), H^1(\partial B) \right)}$ and
$\Vert S^{-1}\Vert_{\mathcal{L}\left(H^1(\partial D), L^2(\partial D) \right)}\leq\epsilon^{-1} \Vert {S}^{-1}\Vert_{\mathcal{L}\left(H^1(\partial B), L^2(\partial B) \right)}$, see 
\cite[Lemma 2.4 and Lemma 2.5]{D-M:MMS:2014} for the last two inequalities.
Using these properties, we deduce that
\begin{equation}\label{key-observation}
2 \pi \Vert S\Vert^{-1}_{\mathcal{L}\left(L^2(\partial B), H^1(\partial B) \right)} \frac{R^2_i(B)}{R_e(B)} R_i(D) \leq C^a_D \leq
4 \pi\Vert S^{-1}\Vert_{\mathcal{L}\left(H^1(\partial B), L^2(\partial B) \right)} \frac{R^2_e(B)}{R_i(B)} R_e(D).
\end{equation}


We can estimate $2 \pi \Vert S\Vert^{-1}_{\mathcal{L}\left(L^2(\partial B), H^1(\partial B) \right)} \frac{R^2_i(B)}{R_e(B)}$ from below by a constant $c(Lip)$ depending only 
the Lipschitz character of $B$ and $4 \pi\Vert S^{-1}\Vert_{\mathcal{L}\left(H^1(\partial B), L^2(\partial B) \right)} \frac{R^2_e(B)}{R_i(B)}$ by a constant $C(Lip)$ also 
depending only on the Lipschitz character of $B$. With these apriori bounds (\ref{key-observation}) becomes
\begin{equation}\label{final-double-estimate}
 c(Lip) R_i(D) \leq C^a_D \leq C(Lip) R_e(D).
\end{equation}
Using Lemma \ref{capacitance-eig-single}, we deduce the following result.
\begin{theorem}\label{convex-obstacles}
Assume that the reference obstacles $B_m,m=1, ...,M$, are convex and satisfy the properties $B(0, R_i(B_m))\subset B_m$. Then,
 there exist two constants $c(Lip)$ and $C(Lip)$ depending only on the Lipschitz character of the obstacles $B_m, m=1,...,M$, such that we have the estimates
 \begin{equation}\label{thikness}
  R_i(D_m) \leq c^{-1}(Lip) (\lambda + 2\mu)^{-1}\lambda^{max}_{eig_m},
 \end{equation}
\begin{equation}\label{lenght}
  R_e(D_m) \geq C^{-1}(Lip) \mu^{-1}\lambda^{min}_{eig_m}.
 \end{equation}
\end{theorem}
The estimate (\ref{thikness})  means that the largest ball contained in $D_m$ has a radius not exceeding $c^{-1}(Lip) (\lambda + 2\mu)^{-1}\lambda^{max}_{eig_m}$. 
Hence (\ref{thikness}) measures the thickness of $D_m$. The estimate (\ref{lenght}) means that the radius of the smallest ball containing $D_m$ is not lower than 
$C^{-1}(Lip) \mu^{-1}\lambda^{min}_{eig_m}$. Hence (\ref{lenght}) measures the length of the obstacle $D_m$.

\begin{center} {\bf{Conclusion}} \end{center}

Based on the asymptotic expansion of the elastic waves by small rigid obstacles, derived in \cite{DPC-SM13-2},
we have shown that any of the two elastic waves, i.e. P-waves or S-waves (precisely SH-waves or SV-waves in 3D elasticity), is enough to localize the obstacles and 
estimate their respective sizes. Compared to the existing literature, see for instance \cite{A-K:2007}, we allow the obstacles to be close and the 
cluster to be spread in any given region. In addition, the derived precise size estimates seem to be new compared to the related literature. 
We stated the MUSIC algorithm based on the mentioned measurements and we believe that performing this algorithm will provide us with accurate numerical tests,
see our previous work \cite{DPC-SM:InvPrbm2012} on point-like obstacles.
Let us also mention that, using our techniques, we can write down a topological derivative based imaging approach from
only the shear or compressional part of the elastic wave, compare to \cite{A-B-G-K-L-W:2013}. We think that these two points deserve to be studied.

\bibliographystyle{abbrv}

\end{document}